\documentclass{ifacconf}

\makeatletter
\let\old@ssect\@ssect 
\makeatother

\usepackage{algpseudocode}
\usepackage{algorithm}
\usepackage{amsmath,amssymb}
\usepackage{graphicx}      
\usepackage{epstopdf}
\usepackage{natbib}        

\newtheorem{lemma}{Lemma}[section]
\newtheorem{corollary}{Corollary}[section]
\newtheorem{remark}{Remark}[section]
\newtheorem{definition}{Definition}[section]

\usepackage{color, hyperref}
\definecolor{darkblue}{rgb}{0.0,0.0,0.6}
\hypersetup{colorlinks,breaklinks,linkcolor=darkblue,urlcolor=darkblue,anchorcolor=darkblue,citecolor=darkblue}
\makeatletter
\def\@ssect#1#2#3#4#5#6{%
	\NR@gettitle{#6}
	\old@ssect{#1}{#2}{#3}{#4}{#5}{#6}
}
\makeatother
\begin{document}

\begin{frontmatter}
	
	\title{On Data-Driven Drawdown Control \\ with Restart Mechanism in Trading} 
	
	\thanks[footnoteinfo]{This paper was supported in part by Ministry of Science and Technology, R.O.C. Taiwan, under Grant: MOST 111--2221--E--007--124--	.}
	

	\author[First]{Chung-Han Hsieh} 
	
	\address[First]{Department of Quantitative Finance, \\National Tsing Hua University, Hsinchu 30044, Taiwan R.O.C. \\(e-mail: \href{ch.hsieh@mx.nthu.edu.tw}{ch.hsieh@mx.nthu.edu.tw}). }


	\begin{abstract}
	This paper extends the existing drawdown modulation control policy to include a novel \textit{restart} mechanism for trading. 
	It is known that the \textit{drawdown modulation policy} guarantees the maximum percentage drawdown no larger than a prespecified drawdown limit for all time with probability one. 
	However, when the prespecified limit is approaching in practice, such a modulation policy becomes a stop-loss order, which may miss  the profitable follow-up  opportunities, if any. 
	Motivated by this,
	we add a data-driven restart mechanism into the drawdown modulation trading system to auto-tune the performance.
	We find that with the restart mechanism, our policy may achieve a superior trading performance to that without the restart, even with a nonzero transaction costs setting. 
	 To support our findings, some  empirical studies using  equity ETF and cryptocurrency with historical price data are provided.
	\end{abstract}

\end{frontmatter}

\section{Introduction} \label{section: introduction}
Starting from the pioneering work by \cite{Markowitz_1952,markowitz1959portfolio}, the portfolio optimization problem is often solved by the \textit{mean-variance} approach. 
That is, the trader seeks an optimal trade-off between payoff and risk measured by the portfolio returns' mean and variance. 
While the variance is widely used as a standard \textit{risk} metric in finance, it is more to the \textit{dispersion risk}, which treats both positive and negative deviation from the mean as equally risky; see \cite{fabozzi2007robustbook} for a good introduction to this topic.

\subsection{Downside Risks}
To remedy the equal riskiness on the dispersion risks, many surrogate risk measures are proposed to pay attention to the \textit{downside risks}.
This includes Value at Risk (VaR); see \cite{jorion2000value}, Conditional Value at Risk (CVaR); see \cite{rockafellar2000optimization}, absolute drawdown; see \cite{magdon2004maximum,hayes2006maximum}, conditional expected drawdown (CED), the tail mean of maximum drawdown distributions, see \cite{goldberg2017drawdown}, and the more general \textit{coherent risk} that axiomatize the risk measures; see \cite{Luenberger_2011, shapiro2021lectures}. 
See also \cite{korn2022drawdown} for an empirical study of comparing various drawdown-based risk metrics. 
In this paper, we focus on a more practical drawdown measure, the \textit{maximum percentage drawdown}, the maximum percentage drops in wealth over time, as the risk measure.

\subsection{Control of Drawdown}
Different types of drawdown and methodologies are studied extensively in the existing literature. 
For example,  optimal drawdown control problem in a continuous-time setting are studied in \cite{grossman1993optimal, cvitanic1994portfolio, chekhlov2005drawdown, Malekpour_Barmish_2012, Malekpour_Barmish_2013}.
See also \cite{boyd2017multi} for a  study on multiperiod portfolio optimization problems involving drawdown as a  constraint in a discrete-time setting. 
A recent study that uses deep reinforcement learning to address practical drawdown issues can be found in \cite{wu2022embedded}.

Among all the existing papers, the prior works in \cite{hsieh2017drawdown,hsieh2017inefficiency} are the most closely related papers, where the key result, the so-called \textit{drawdown modulation lemma}, is proved. 
Roughly speaking, it shows a necessary and sufficient conditions for a broad class of control policies that guarantees an almost sure maximum percentage drawdown protection for all time. 
However, \cite{hsieh2017drawdown} indicate that, in practice, when the prespecified drawdown limit is approaching, the trading policy may behave like a stop-loss order, and the trading may be stopped; see also \cite{hsieh2022generalization} for a study on a class of affine policies with a stop-loss order.
To this end, we extend the drawdown modulation policy with a novel restart mechanism to remedy the stop-loss phenomenon.

\subsection{Contributions of the Paper}
This paper extends the existing drawdown modulation policy with a novel restart mechanism.
The preliminaries are provided in Section~\ref{section: preliminiaries}.
We formulate an optimal drawdown control problem for a two-asset portfolio with one risky and one riskless asset.
Then we extend the existing drawdown modulation theory so that the riskless asset is explicitly involved; see Lemma~\ref{lemma: drawdown modulation}. The necessary and sufficient conditions for a broad class of control policy which we call the drawdown modulation policy are provided.
Then, the modulation policy with a restart mechanism is discussed in Section~\ref{section: drawdown modulaiton with restart}.
The idea of the restart is simple: When the percentage drawdown up-to-date is close to the prespecified drawdown limit, we restart the trades with an updated policy.
We also provide numerical examples involving ETF and cryptocurrency historical price data to support our findings; see Section~\ref{section: empirical studies}.

\medskip
\section{Preliminaries} \label{section: preliminiaries}
We now provide some useful preliminaries for the sections to follow.

\subsection{Asset Trading Formulation}
Fix an integer $N>0.$
For stage $k=0, 1,  \dots, N$, we let~$S(k)>0$ denote the prices of the underlying financial asset at stage $k$. The associated \textit{per-period returns} are given by
$
X(k):= \frac{S(k+1) - S(k)}{S(k)}
$
and the returns are assumed to be  bounded, i.e.,
$
X_{\min} \leq X(k) \leq X_{\max} 
$
with $X_{\min}$ and $X_{\max}$ being points in the support, denoted by $\mathcal{X}$,
and satisfying
$
-1 < X_{\min} < 0 < X_{\max}.
$
For the money market asset, e.g., bond or bank account, we use $r_f(k)$ to denote the interest rate at stage~$k$.

\medskip
\begin{remark} \rm
 Note that the returns considered here are not necessarily independent and can have an arbitrary but  bounded distribution with the bounds $X_{\min}$ and $X_{\max}$.
\end{remark}

\subsection{Account Value Dynamics} 
Beginning at some initial account value $V(0) > 0$, consider a portfolio consisting of two assets, with one being risky and the other being a riskless asset with interest rate~$r_f(k) \in [0, X_{\max}]$ for all $k$ almost surely. 
For stage~$k=0,1, \dots$, we let $V(k)$ denote the account value at stage~$k$. Then the evolution of the account value dynamics is described by the stochastic recursion
\begin{align*}
	V(k+1)=V(k)+u(k)X(k) + (V(k) - u(k))r_f(k) .
\end{align*}
Given a \textit{prespecified drawdown limit} $d_{\max} \in (0,1)$, we focus on conditions on selecting a policy $u(k)$ under which satisfaction of the constraint
$
d(k)\leq d_{\max}
$
is assured for all~$k$ with probability one where $d(k)$ is the percentage drawdown up to date $k$, which is defined below.

\medskip
\begin{definition}[Maximum Percentage Drawdown]
	For $k=0,1, \dots, N$, 
	the \textit{percentage drawdown} up to date $k$, denoted by $d(k)$, is defined as
	$$
	d(k) := \frac{V_{\max}(k) - V(k)}{V_{\max}(k)}
	$$
	where
	$
	V_{\max}(k) := \max_{0\leq i\leq k} V(i).
	$
	The \textit{maximum percentage drawdown}, call it $d^*$, is then defined as
	$$
	d^* := \max_{0\leq k\leq N} d(k).
	$$
\end{definition}

\medskip
\begin{remark} \rm
	It is readily seen that the percentage drawdown satisfies $d(k)\in [0, 1]$ for all $k$ with probability one.
\end{remark}

\medskip
\section{Drawdown Modulation with Restart} \label{section: drawdown modulaiton with restart}
According to \cite{hsieh2017drawdown}, it states a necessary and sufficient condition on any trading policy~$u(k)$ that guarantees the percentage drawdown up to date $d(k)$ is no greater than a given level $d_{\max}$ for all $k$ with probability one. Below, we extend the result to include a riskless asset.

\medskip
\begin{lemma}[Drawdown Modulation] \label{lemma: drawdown modulation}
{\it 	Let $d_{\max} \in (0,1)$ be given.
 An trading policy $u(\cdot)$ guarantees prespecified drawdown limit satisfying $d(k) \leq d_{\max}$ for all $k$ with probability one if and only if for all $k$, the condition
$$
-\frac{ M(k) + r_f(k) }{X_{\max} - r_f(k)}V(k)\leq u(k) \leq \frac{M(k) + r_f(k)}{|X_{\min}| + r_f(k)  }V(k)
$$
is satisfied along all sample paths
where 
$$
M(k) := \frac{d_{\max} - d(k)}{1-d(k)}.
$$}
\end{lemma}

\textit{Proof.} 
The idea of the proof is similar to that of \cite{hsieh2017drawdown}. 
However, for the sake of completeness, we provide full proof here.
To prove necessity, assuming that
	$d(k) \le d_{\max}$  for all~$k$ and all sequences of returns, we must show the condition on~$u(k)$ holds along all sequences of returns. 
	Fix~$k$. Since both~$d(k) \le d_{\max}$ and $d(k+1) \leq d_{\max}$  for all sequences of returns, we claim this forces the required inequalities on~$u(k)$. 
	Without loss of generality, we prove the right-hand inequality for the case $u(k) \ge 0$ and note that an almost identical proof also works for $u(k) < 0$.  
	To establish the condition on~$u(k)$ for all sequences of returns, it suffices to consider the path with the worst loss~$|X_{\min}|u(k)$. In this case,  we have $V_{\max}(k+1) = V_{\max}(k)$. Hence,
	\begin{align*}
		& d(k + 1) \\
		&= \frac{V_{\max}(k+1) - V(k+1)}{V_{\max}(k+1)}\\
		&= \frac{V_{\max}(k) - V(k+1)}{V_{\max}(k)}\\
		&= \frac{V_{\max}(k) - V(k) + u(k)|X_{\min}| - (V(k) - u(k))r_f(k)}{V_{\max}(k)}\\		
		&= d(k) + \frac{u(k)|X_{\min}| - (V(k) - u(k))r_f(k)}{V_{\max}(k)} \leq d_{\max}
	\end{align*}
	We now substitute
	$
	{V_{\max }}(k) = \frac{{V( k )}}{{1 - d( k )}} >0
	$
	into the inequality above and obtain
	\[
	| X_{\min}  | u(k) - (V(k) - u(k)) r_f(k) \leq M(k)V(k),
	\]
	where $M(k) = \frac{d_{\max} - d(k)}{1- d(k)}$.
	This implies that
	\[
	(	| X_{\min}  | + r_f(k)) u(k) \leq \left( M(k) + r_f(k) \right) V(k).
	\]
	Or equivalently,
	\[
	u(k) \leq  \frac{ M(k) + r_f(k)}{ | X_{\min } | + r_f(k) }V(k).
	\]  
	To prove sufficiency,  assuming that the stated bounds on~$u(k)$ hold along all sequences of returns, 
	we must show~$d(k) \le d_{\max}$  for all~$k$  and all sequences of returns.  Proceeding by induction,  for~$k=0$, we trivially have~{$d(0)=0 \leq d_{\max}$}.  To complete the inductive argument, we assume that~$d(k) \le d_{\max}$  for all sequences of returns, and  must  show~$d(k+1) \le d_{\max}$  for all sequences of returns. Without loss of generality, we again provide a proof for the case~{$u(k) \ge 0$} and note that a nearly identical proof is used for~{$u(k) < 0$}. 
	Indeed, by noting~that
	\begin{align*}
		d( k + 1 ) 
		&= 1 - \frac{{V( k + 1 )}}{{{V_{\max }}( k + 1 )}},
	\end{align*}
	and $V_{\max}(k) \leq V_{\max}(k+1)$  for all sequences of returns, we split the argument into two cases: If $V_{\max}(k) < V_{\max}(k+1)$, then $V_{\max}(k+1) = V(k+1)$ and we have 
	$
	d(k+1)=0 \leq d_{\max}.
	$ 
	On the other hand, if $V_{\max}(k) = V_{\max}(k+1)$, with the aid of the dynamics of account value, we~have
	\begin{align*}
		d( k + 1 ) 
		&= 1 - \frac{ V( k ) + u( k )X(k) + (V(k) - u(k)) r_f(k) }{V_{\max }( k )}\\ 
		&\leq  1 + \frac{ -V( k )(1 + r_f(k))  + u( k )(|X_{\min}| + r_f(k))  }{V_{\max }( k )}.
	\end{align*}
	Using the upper bound on $u(k)$; i.e.,
	\[
		u(k) \leq  \frac{ M(k) + r_f(k)}{ | X_{\min } | + r_f(k) }V(k)
	\] and  $
	{V_{\max }}(k) = \frac{{V( k )}}{{1 - d( k )}}, 
	$
	  we obtain 
	\begin{align*}
	d( k + 1 ) 
	&\leq  1 +  (-1 + M(k)) (1-d(k))   \\
	&=  d(k)  + M(k)(1-d(k))  .
\end{align*}
Using the definition of $M(k) = \frac{d_{\max} - d(k)}{1-d(k)}$,  it follows that~$
	d( k + 1 ) 
	 \leq d_{\max},
$
and the proof is complete. 	 \hfill  \qed

\medskip
\subsection{Drawdown Modulation Policy}
Consistent with \cite{hsieh2017drawdown}, fix the prespecified drawdown limit $d_{\max} \in (0,1)$. 
With the aid of Lemma~\ref{lemma: drawdown modulation}, one can readily obtain a class of policy functions $u(k)$ expressed as a \textit{linear time-varying feedback controller} parameterized by a gain~$\gamma$,  leading to the satisfaction of the drawdown specification. Specifically,
we express $u(k)$ in the feedback form 
\begin{align} \label{eq: modulator }
u(k):= K(k)V(k) 
\end{align}
with $K(k):=\gamma  M(k)$ and  
$$
\gamma \in \Gamma := \left[ \frac{-1}{X_{\max} - \max_k r_f(k) },   \, \frac{1}{ | X_{\min} | + \max_k r_f(k) } \right] .
$$
Equation~(\ref{eq: modulator }) is called the \textit{drawdown modulation policy}, which is parameterized by the two parameters~$(\gamma, d_{\max})$.

\medskip
\begin{remark}
	$(i)$ It is readily verified that the drawdown modulation policy~(\ref{eq: modulator }) satisfies Lemma~\ref{lemma: drawdown modulation}. 
	$(ii)$ To link back to finance concepts, the special case of \textit{buy-and-hold} is obtained when $K(k)  \equiv 1$.
	Note that $u(k) < 0$ stands for \textit{short selling}.\footnote{Short selling a stock means that  a trader borrows the stocks from someone who owns it and selling it with the hope that the prices of the stock will drop in the near future; see \cite{Luenberger_2011}.} 
	$(iii)$ Instead of using a fixed feasible set~$\Gamma$, it is also possible to allow a time-varying feasible set, say~$\Gamma_k$, to reflect the time dependency of the returns.
\end{remark}

\medskip
 \begin{corollary}[Maximum Drawdown Protection] 
	{\it With the drawdown modulation policy~{\rm (\ref{eq: modulator })}, the maximum percentage drawdown satisfies $d^* \leq d_{\max}$.}
\end{corollary}

\textit{Proof.}
	Since the drawdown modulation policy satisfies Lemma~\ref{lemma: drawdown modulation}, it assures $d(k) \leq d_{\max}$ for all $k$ with probability one. Therefore, it follows that
	\[
	d^* = \max_{0\leq k\leq N} d(k) \leq d_{\max}. \qed
	\]

\medskip
\subsection{Optimal Drawdown Control Problem}
Having obtained the drawdown modulation policy~(\ref{eq: modulator }), a natural question arises of how to select an ``optimal"~$\gamma.$ To this end, we define the \textit{total return }up to terminal stage~$N$ as a ratio
\[
R_\gamma (N) := \frac{V(N)}{V(0)}
\]
where the subscript of $R_\gamma(\cdot)$ is used to emphasize the dependence on the gain $\gamma$.
Define $J(\gamma) := 	  \mathbb{E}[ R_\gamma (N) ]$.
Then,  we consider a multiperiod drawdown-based  stochastic optimization problem  
\begin{align*}
J^* :=	\max_{\gamma \in \Gamma} J(\gamma) 
\end{align*}
subject to  
\begin{align*}
V(k+1) 
	&= V(k)+u(k)X(k) + (V(k) - u(k))r_f(k) \\
	& = [1+ r_f(k) + \gamma  M(k) (X(k)  -  r_f(k)) ]V(k).
\end{align*}
It is readily verified that 
$$
\frac{V(N)}{V(0)} =  \prod_{k=0}^{N-1} [1+ r_f(k) + \gamma  M(k) (X(k)  -  r_f(k)) ].
$$
 Therefore, we rewrite the problem as the following equivalent form
\begin{align} \label{problem: drawdown control}
	&\max_{\gamma \in \Gamma} \mathbb{E}[ R_\gamma (N) ] \nonumber \\
	&= 	\max_{\gamma \in \Gamma} \mathbb{E}\left[ \prod_{k=0}^{N-1} [1+ r_f(k) + \gamma M(k) (X(k)  -  r_f(k)) ] \right] .
\end{align}
In the sequel, we shall use $\gamma^*$ to denote a maximizer of the optimization problem above.
In practice, if one view that the optimum~$\gamma^*$ obtained may be too aggressive, a practical but arguably suboptimal way is to introduce an additional \textit{fraction}, call it $\alpha$, that is used to shrink the investment  size; see \cite{maclean2010long} for a similar idea.
That is, instead of working with~$\gamma^*$, one may work with $\alpha \gamma^*$ where~$\alpha \in (0, 1]$. 

\medskip
\begin{remark}[Non-Convexity] \rm
 It is important to note that solving Problem~(\ref{problem: drawdown control}) is challenging since the modulation function $M(k)$ depends on $\gamma$ and the history of~$X(0), \dots, X(k-1)$, which in general yields a nonconvex problem; e.g., see  Figure~\ref{fig:jgammainsample} in Section~\ref{section: empirical studies} for an illustration of non-convexity nature.
  Therefore, Monte-Carlo simulations are often needed to obtain the optimum. 
\end{remark}
 
\medskip
\subsection{Modulation with Restart}
As mentioned in Section~\ref{section: introduction}, while the derived drawdown modulation policy  provides almost sure drawdown protection, it may incur a stop-loss behavior.
To remedy this, we now introduce a restart mechanism into the modulation policy.
 Specifically, let $\varepsilon \in (0,d_{\max})$ be a prespecified threshold parameter. Then we set the \textit{threshold} for restarting the trade~by
\begin{align} \label{ineq: restart threshold}
	 d(k) + \varepsilon > d_{\max}
\end{align}
If, at some stage $k = k_0$, Inequality~(\ref{ineq: restart threshold}) is satisfied, the trading is restarted by re-initializing $d(k_0):=0$ and reset the time-varying gain function $K(k)$ of the modulation policy~$u(k) =K(k)V(k)$  at that stage $k=k_0$ as 
\begin{align}\label{eq: time-varying restart gain}
	K(k_0) := \gamma^* \alpha e^{-k_0/N} M(k_0)
\end{align}
where $ \alpha e^{-k_0/N}$ represents a \textit{forgetting factor} with fraction~$\alpha \in (0, 1]$ mentioned previously.
  Then we continue the trade until the next restart stage or to the terminal stage $N$.

\medskip
\begin{remark} \rm
$(i)$ The forgetting factor in Equation~(\ref{eq: time-varying restart gain}) reflects the idea that  the trading size should be shrunk after the restart. 
Said another way, if the trades are approaching the prespecified drawdown limit $d_{\max}$,  the follow-up trades should be more conservative after the restart. 
$(ii)$ Note that after the restart, the control policy satisfies $|u(k_0)| \leq  |\gamma^*| \alpha d_{\max} $ since $M(k) \leq d_{\max}$ for all $k$ with probability one.
$(iii)$ While it does not consider in this paper, we should note that the optimal $\gamma^*$ in Equation~(\ref{eq: time-varying restart gain}) can also be re-calculated at each restart time by using the previous~$k_0-m$ stages information for some integer $m > 1$; see also \cite{wang2022data} for a similar idea for obtaining a data-driven log-optimal portfolio via a sliding-window approach.
\end{remark}

\medskip
\section{Illustrative Examples} \label{section: empirical studies}
In this section, two trading examples are proposed to support our theory. The first example is trading with ETF and riskless asset. 
The second example is trading with Bitcoin and  a riskless asset. 
For the sake of simplicity,
we take a constant daily interest rate $r_f(k) := 0.01/365$ for all $k$, which corresponds to a  $1\%$ annual rate.
While our theory allows \textit{leveraging}, in the sequel, we impose an additional \textit{cash-financing} condition by imposing that~$|u(k)| \leq V(k)$ for all $k$, which corresponds to~$|K^*(k)| \leq 1.$

\medskip
\subsection{Trading with ETF and Riskless Asset} \label{subsection: trading with ETF}
\begin{figure}[h!]
	\centering
	\includegraphics[width=0.9\linewidth]{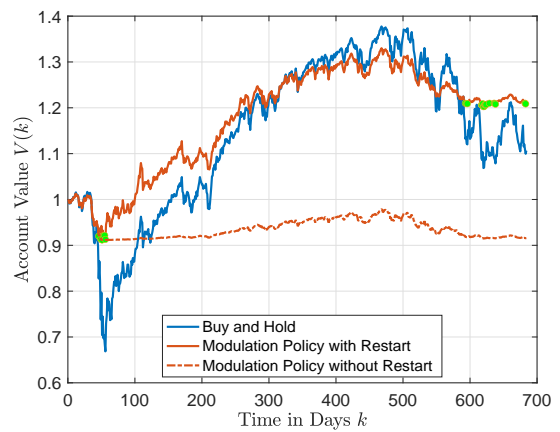}
	\caption{Drawdown Modulation with/without Restart (Green Dots Indicate a Restart)}
	\label{fig:vttradingresults}
\end{figure}

Consider the Vanguard Total World Stock Index Fund ETF (Ticker: VT)\footnote{VT invests in both foreign and U.S. stocks; hence it can be viewed as a good representative of the global stock market} to be the risky asset covering a one-year in-sample period for optimization from~January~02,~2019 to January 02, 2020, and the out-of-sample period from~January 02, 2020 to September 20, 2022, which contains a total of $N=684$ trading days; see Figure~\ref{fig:vttradingresults} where the cyan colored trajectory is used for in-sample optimization and the blue colored trajectory is used for the out-of-sample trading test.  
It should be noted that due to the  COVID-19 pandemic, the considered prices covering the 2020 stock market crash from  February~20,~2020 to  April 7, 2020, and a recovery period after the crash. 
Thus, we view this dataset as an excellent back-testing case for our proposed drawdown modulation policy with the~restart.

\begin{figure}[h!]
	\centering
	\includegraphics[width=0.9\linewidth]{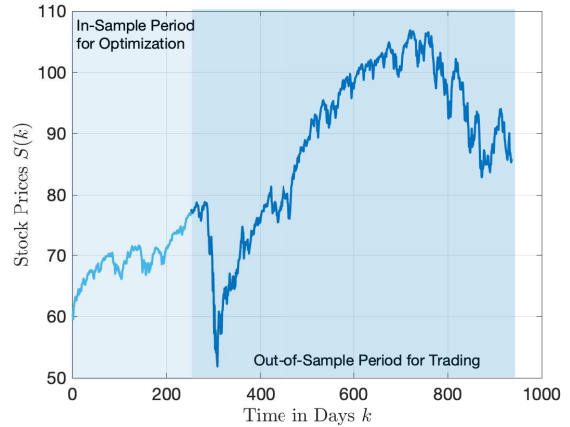}
	\caption{Stock Prices of VT}
	\label{fig:vtprices20202022}
\end{figure}

Without loss of generality, consider the initial account value to be $V(0) := \$1$. To implement the drawdown modulation policy with restart, we set $d_{\max}:=0.1$ and restart threshold $\varepsilon:= d_{\max}/10$.
With the data collected in the in-sample period, the corresponding feasible set is~$\Gamma \approx [-30.79, \; 34.9]$. 
Then we numerically solve the optimization problem~(\ref{problem: drawdown control})  via Monte-Carlo simulations.
It follows that any $\gamma \in (8.5, 34.9)$ share an almost identical optimal value when the cash-financing condition is imposed. 
For the sake of risk-averseness, we pick the infimum of the candidates, i.e.,~$\gamma^* := 8.5$.
Using this $\gamma$, we obtain the drawdown modulation policy $u^*(k) = \gamma^* M(k) V(k).$

The corresponding trading performance is shown in Figure~\ref{fig:vttradingresults}, where the green dots indicate that the trade is restarted. 
In the same figure, we also compare it with the standard buy-and-hold strategy\footnote{Here, we mean buy and hold on the risky asset VT with $K(k)=1$.} and the modulation policy without restart. 
We see clearly that the modulation policy with restart leads to superior performance.
Some  key performance metrics, including maximum drawdown and cumulative rate returns, and $N$-period Sharpe ratio\footnote{The per-period Sharpe ratio is $SR:= \frac{ \mu - r_{f}}{\sigma}$ where $\mu $ is the per-period sample mean return,  $\sigma$ is the per-period sample standard deviation, and $r_{f}$ is the per-period riskless return.}   are reported in Table~\ref{table:Performance Metrics for VT}.

%

\begin{table*}[h!]
	\centering
	\caption{Performance Metrics for Trading VT and Riskless Asset}
	\begin{tabular}{ | c | c | c | c |} 
		\hline
		\multicolumn{4}{c}{Trading Performance with $d_{\max}=0.1$ and $\varepsilon = d_{\max}/10$}  \\ 
		\hline\hline
		Trading Policy  & Buy and Hold & Modulation without Restart & Modulation with Restart  \\ 
		\hline
		Maximum percentage drawdown $d^*$ & $34.24\%$ & $10.00\%$ & $\textbf{10.00}\%$ \\ 
		\hline
		Cumulative rate of return $\frac{V(N) - V(0)}{V(0)}$ & $10.64\%$ & $-8.422\%$ & $\textbf{20.97}\%$   \\ 
		\hline
		$N$-period Sharpe ratio $\sqrt{N}\cdot SR$ & $ 0.4082$ & $-1.3375$ & $\textbf{1.0462}$  \\ 
		\hline
	\end{tabular}
	\label{table:Performance Metrics for VT}
\end{table*}

\medskip
\subsection{Trading with Cryptocurrency and Riskless Asset}
\begin{figure}[h!]
	\centering
	\includegraphics[width=0.9\linewidth]{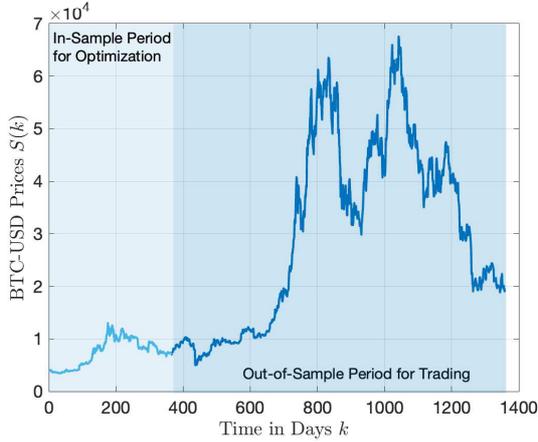}
	\caption{Prices of BTC-USD (In-Sample and Out-of-Sample)}
	\label{fig:btcusdinandoutsampleprices}
\end{figure}

As a second example, we consider a portfolio consisting of cryptocurrency BTC-USD and a riskless asset. The BTC-USD asset covers the same in-sample and out-of-sample periods described in Example~\ref{subsection: trading with ETF}. 
From January~02,~2020 to September 20, 2022, it has a total of~$N=993$ trading days\footnote{It is worth noting that  the cryptocurrency is typically traded at~24 hours a day, seven days a week. 
	Therefore, it has longer tradings days than that trades with VT in Section~\ref{subsection: trading with ETF}.}; see Figure~\ref{fig:btcusdinandoutsampleprices}. The corresponding feasible set for~$\gamma$ is $\Gamma = [-5.761, \; 7.083]$.

\begin{figure}[h!]
	\centering
	\includegraphics[width=0.9\linewidth]{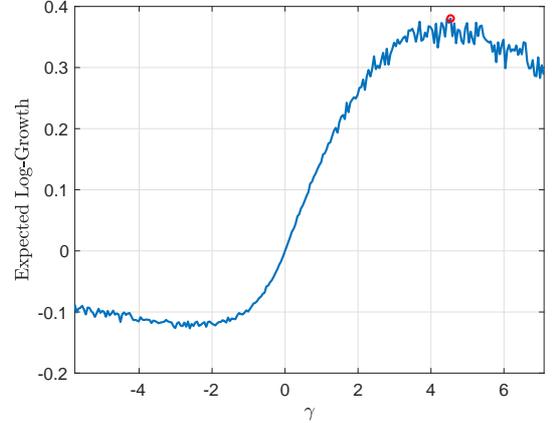}
	\caption{Seeking the Optimum: $J(\gamma) $ Versus $\gamma \in \Gamma$}
	\label{fig:jgammainsample}
\end{figure}

\begin{figure}[h!]
	\centering
	\includegraphics[width=0.9\linewidth]{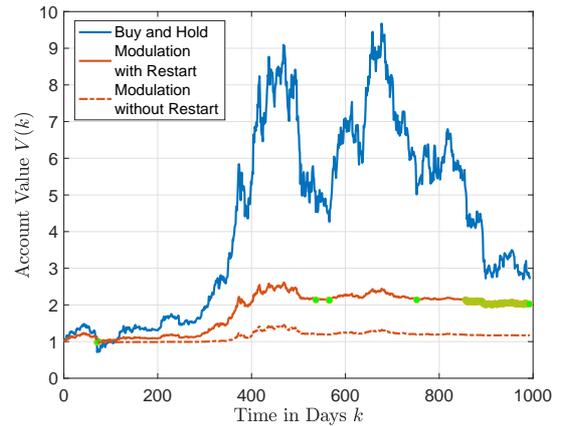}
	\caption{Trading Performance of BTC-USD with $d_{\max}:=0.2$ and $\varepsilon:= d_{\max}/10$.}
	\label{fig:btcusdtradingex}
\end{figure}

Take $d_{\max} := 0.2$ and restart threshold $\varepsilon:=d_{\max}/10$.
By solving the optimal drawdown control problem~(\ref{problem: drawdown control}), we obtain $\gamma^* \approx 5.138$; see  Figure~\ref{fig:jgammainsample} for $J(\gamma)$ versus $\gamma \in \Gamma$. 
Note that the~$J(\gamma)$ is clearly not concave for~$\gamma \in \Gamma$.
To consider the volatile nature of cryptocurrency and  unforeseen estimation error, we consider using a fractional~$\gamma^*$ by~$\alpha \gamma^*$ with $\alpha = 1/2$. 
That is, $ u^*(k) = K^*(k) V(k)$ with $K^*(k) = \frac{\gamma^*}{2} M(k) V(k)$.
The trading performances using drawdown modulation policy with and without restart, and buy-and-hold strategy are shown in Figure~\ref{fig:btcusdtradingex}, where the green dots indicate that the trades were restarted. Some performance metrics are summarized in Table~\ref{table:Performance Metrics for BTC-USD} where we see that the modulation with restart provides the highest Sharpe ratio among all the other strategies.

\begin{table*}[h!]
	\centering
	\caption{Performance Metrics for Trading BTC-USD and Riskless Asset}
	\begin{tabular}{ | c | c | c | c |} 
		\hline
		\multicolumn{4}{c}{Trading Performance with $d_{\max} = 0.2$ and $\varepsilon = d_{\max}/10$ }  \\ 
		\hline\hline
		Trading Policy  & Buy and Hold & Modulation without Restart & Modulation with Restart  \\ 
		\hline
		Maximum percentage drawdown $d^*$ & $72.12\%$ & $20.00\%$ & $\textbf{20.18}\%$ \\ 
		\hline
		Cumulative rate of return $\frac{V(N) - V(0)}{V(0)}$ & $170.43\%$ & $17.41\%$ & $\textbf{101.82}\%$   \\ 
		\hline
		$N$-period Sharpe ratio $\sqrt{N}\cdot SR$ & $1.4424$ & $0.6970$ & $\textbf{2.1363}$  \\ 
		\hline
	\end{tabular}
	\label{table:Performance Metrics for BTC-USD}
\end{table*}

\subsubsection{Transaction Costs Considerations.}
\begin{figure}[h!]
	\centering
	\includegraphics[width=0.9\linewidth]{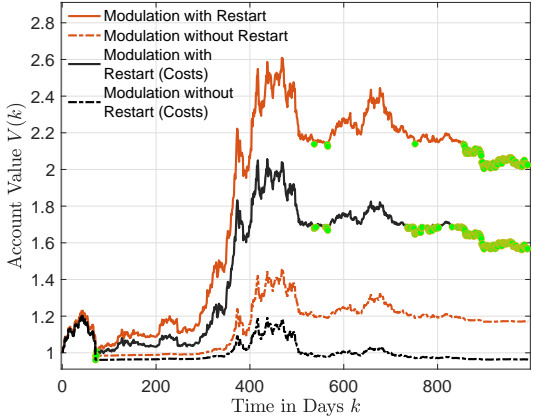}
	\caption{Trading Performance with and without Transaction Costs}
	\label{fig:btcusdtradingexcostscomparison}
\end{figure}

Traditionally, with cryptocurrency as the underlying asset, most exchanges typically charge transaction fees from $0\%$ to $1\%$, depending on the trading volumes. 
While we are not developing the drawdown modulation theory involving the transaction costs, in this example,  we  perform backtests to see the effects with percentage costs of $0.1\%$ per trade.\footnote{ 
This cost level is typical in common cryptocurrency exchanges such as Binance; see \cite{hsieh2022robustness}.}
Take~$d_{\max} = 0.2$, $\varepsilon = d_{\max}/10$ and $\gamma^* = 2.69$, 
the trading performances are shown in Figure~\ref{fig:btcusdtradingexcostscomparison}. 
From the figure, we see that, even with transaction costs, the proposed modulation policy with a restart mechanism is still superior to that without a restart. The corresponding Sharpe ratio is about $1.521$ compared with that without restart, where the Sharpe ratio is down to $-0.043$.

\medskip
\section{Concluding Remarks}
This paper extended the existing drawdown modulation policy with a novel data-driven restart mechanism. 
We first derived a drawdown modulation policy by proving a lemma for a two-asset portfolio involving one risky and riskless asset.
Then we modified the derived modulation policy by adding the restart mechanism.
Subsequently, we showed an auto-tuned trading performance using historical price data for an ETF and cryptocurrency within a duration where the prices are volatile. Overall, our modulation policy with restart dynamically controls the drawdown and maintains the profits level.

Some possible future research directions might include the multi-asset problem, which will be considered in the subsequent  journal version of this paper.
The other possibility is to consider the consumption-investment problem; i.e., 
$$
\sup_\gamma \mathbb{E}\left[ \sum_{k=0}^{T-1} U_1(c(k)) + U_2(V(N)) \right]
$$
 where $U_1$ and $U_2$ are two utility functions with $c(k)$ being the amount consumed at time $k$. See \cite{cvitanic2004introduction}.
Lastly,  one might consider  a more general portfolio optimization problem that seeks $\gamma$ maximizing~$
 \sup_\gamma \mathbb{E}\left[ U(V_\gamma (N)) \right]
$
where $U: \mathbb{R} \to \mathbb{R}$ is a strictly increasing and concave utility function. One example is to consider a general \textit{hyperbolic absolute risk aversion}~(HARA) class of utilities;~i.e.,
$$
U(x) = \frac{1 -\theta}{\theta} \left( \frac{ax}{1-\theta} + b \right)^\theta 
$$ 
for~$a,b>0$ and $\theta < 1$.
The HARA utility functions include the quadratic, exponential, and the isoelastic utility function as  special cases.

\medskip
\section*{Acknowledgment}
The author thanks Chia-Yin Lee for coding and running some preliminary numerical examples on the early draft of this paper.

\appendix
\bibliography{refs}

\end{document}